\documentclass[12pt]{article}
\usepackage{bbm}
\usepackage{amsfonts}
\usepackage{mathrsfs}
\usepackage{amscd}
\usepackage{amsmath,amsfonts,amssymb,amscd}
\usepackage{indentfirst,graphics,epsfig,psfrag}
\input{epsf}
\usepackage{ifpdf}
\usepackage{enumerate}
\usepackage{appendix}
\usepackage{enumerate}

\usepackage{lineno}

\makeatletter \@addtoreset{figure}{section} \makeatother
\makeatletter
\long\def\@makecaption#1#2{%
   \vskip 10\p@
   \setbox\@tempboxa\hbox{{#1}\ \ #2}%
   \ifdim \wd\@tempboxa >\hsize
       {#1}\ \ #2\par
   \else
       \hbox to\hsize{\hfil\box\@tempboxa\hfil}%
   \fi}
\makeatother

\newtheorem{thm}{Theorem}[section]
\newtheorem{cor}[thm]{Corollary}
\newtheorem{lem}[thm]{Lemma}

\newtheorem{con}[thm]{Conjecture}

\newtheorem{obs}[thm]{Observation}
\newtheorem{pro}[thm]{Proposition}
\newtheorem{prob}[thm]{Problem}
\newtheorem{que}[thm]{Question}

\newcommand{\qed}{{\hfill\rule{3pt}{7pt}}}

\def\qed{\hfill \rule{4pt}{7pt}}

\begin{document}
\title{Rainbow connections of graphs -- A survey\footnote{Supported by NSFC.}}
\author{
\small  Xueliang Li, Yuefang Sun\\
\small Center for Combinatorics and LPMC-TJKLC\\
\small Nankai University, Tianjin 300071, P.R. China\\
\small E-mails: lxl@nankai.edu.cn, syf@cfc.nankai.edu.cn
 }
\date{}
\maketitle
\begin{abstract}
The concept of rainbow connection was introduced by Chartrand et al.
in 2008. It is fairly interesting and recently quite a lot papers
have been published about it. In this survey we attempt to bring
together most of the results and papers that dealt with it. We begin
with an introduction, and then try to organize the work into five
categories, including (strong) rainbow connection number, rainbow
$k$-connectivity, $k$-rainbow index, rainbow vertex-connection
number, algorithms and computational complexity. This survey also
contains some conjectures, open problems or questions.\\[2mm]
\end{abstract}

\section{Introduction}

\subsection{Motivation and definitions}
Connectivity is perhaps the most fundamental graph-theoretic
subject, both in combinatorial sense and the algorithmic sense.
There are many elegant and powerful results on connectivity in graph
theory. There are also many ways to strengthen the connectivity
concept, such as requiring hamiltonicity, $k$-connectivity, imposing
bounds on the diameter, and so on. An interesting way to strengthen
the connectivity requirement, the rainbow connection, was introduced
by Chartrand, Johns, McKeon and Zhang \cite{Chartrand 1} in 2008,
which is restated as follows:

This new concept comes from the communication of information between
agencies of government. The Department of Homeland Security of USA
was created in 2003 in response to the weaknesses discovered in the
transfer of classified information after the September 11, 2001
terrorist attacks. Ericksen \cite{Ericksen} made the following
observation: An unanticipated aftermath of those deadly attacks was~
the realization that law enforcement and intelligence agencies
couldn't communicate with each other through their regular channels,
from radio systems to databases. The technologies utilized were
separate entities and prohibited shared access, meaning that there
was no way for officers and agents to cross check information
between various organizations.

While the information needs to be protected since it relates to
national security, there must also be procedures that permit access
between appropriate parties. This two-fold issue can be addressed by
assigning information transfer paths between agencies which may have
other agencies as intermediaries while requiring a large enough
number of passwords and firewalls that is prohibitive to intruders,
yet small enough to manage (that is, enough so that one or more
paths between every pair of agencies have no password repeated). An
immediate question arises: What is the minimum number of passwords
or firewalls needed that allows one or more secure paths between
every two agencies so that the passwords along each path are
distinct?

This situation can be modeled by graph-theoretic model. Let $G$ be a
nontrivial connected graph on which an edge-coloring $c:
E(G)\rightarrow \{1,2,\cdots,n\}$, $n\in \mathbb{N}$, is defined,
where adjacent edges may be colored the same. A path is $rainbow$ if
no two edges of it are colored the same. An edge-coloring graph $G$
is $rainbow ~connected$ if any two vertices are connected by a
rainbow path. An edge-coloring under which $G$ is rainbow connected
is called a $rainbow$ $coloring$. Clearly, if a graph is rainbow
connected, it must be connected. Conversely, any connected graph has
a trivial edge-coloring that makes it rainbow connected; just color
each edge with a distinct color. Thus, we define the
$rainbow~connection~number$ of a connected graph $G$, denoted by
$rc(G)$, as the smallest number of colors that are needed in order
to make $G$ rainbow connected \cite{Chartrand 1}. A rainbow coloring
using $rc(G)$ colors is called a $minimum$ $rainbow$ $coloring$. So
the question mentioned above can be modeled by means of computing
the value of rainbow connection number. By definition, if $H$ is a
connected spanning subgraph of $G$, then $rc(G)\leq rc(H)$. For a
basic introduction to the topic, we refer the readers to Chapter 11
in \cite{Chartrand 5}.

In addition to regarding as a natural combinatorial measure and its
application for the secure transfer of classified information
between agencies, rainbow connection number can also be motivated by
its interesting interpretation in the area of networking
\cite{Chakraborty}: Suppose that $G$ represents a network (e.g., a
cellular network). We wish to route messages between any two
vertices in a pipeline, and require that each link on the route
between the vertices (namely, each edge on the path) is assigned a
distinct channel (e.g. a distinct frequency). Clearly, we want to
minimize the number of distinct channels that we use in our network.
This number is precisely $rc(G)$.

Let $c$ be a rainbow coloring of a connected graph $G$. For any two
vertices $u$ and $v$ of $G$, a $rainbow$~$u-v$~$geodesic$ in $G$ is
a rainbow $u-v$ path of length $d(u,v)$, where $d(u,v)$ is the
distance between $u$ and $v$ in $G$. A graph $G$ is
$strong$~$rainbow$~$connected$ if there exists a rainbow $u-v$
geodesic for any two vertices $u$ and $v$ in $G$. In this case, the
coloring $c$ is called a $strong$~$rainbow$~$coloring$ of $G$.
Similarly, we define the $strong~rainbow~connection~number$ of a
connected graph $G$, denoted $src(G)$, as the smallest number of
colors that are needed in order to make $G$ strong rainbow connected
\cite{Chartrand 1}. Note that this number is also called the
$rainbow$ $diameter$ $number$ in \cite{Chakraborty}. A strong
rainbow coloring of $G$ using $src(G)$ colors is called a
$minimum~strong~rainbow~coloring$ of $G$. Clearly, we have
$diam(G)\leq rc(G)\leq src(G)\leq m$, where $diam(G)$ denotes the
diameter of $G$ and $m$ is the size of $G$.

In a rainbow coloring, we only need to find one rainbow path
connecting any two vertices. So there is a natural generalizaiton:
the number of rainbow paths between any two vertices is at least an
integer $k$ with $k\geq 1$ in some edge-coloring. A well-known
theorem of Whitney \cite{Whitney} shows that in every
$\kappa$-connected graph $G$ with $\kappa \geq 1$, there are $k$
internally disjoint $u-v$ paths connecting any two distinct vertices
$u$ and $v$ for every integer $k$ with $1\leq k\leq \kappa$. Similar
to rainbow coloring, we call an edge-coloring a $rainbow$
$k$-$coloring$ if there are at least $k$ internally disjoint $u-v$
paths connecting any two distinct vertices $u$ and $v$. Chartrand,
Johns, McKeon and Zhang \cite{Chartrand 2} defined the $rainbow$~
$k$-$connectivity$ $rc_k(G)$ of $G$ to be the minimum integer $j$
such that there exists a $j$-edge-coloring which is a $rainbow$
$k$-$coloring$. A rainbow $k$-coloring using $rc_k(G)$ colors is
called a $minimum$ $rainbow$ $k$-$coloring$. By definition,
$rc_k(G)$ is the generalization of $rc(G)$ and $rc_1(G)=rc(G)$ is
the rainbow connection number of $G$. By coloring the edges of $G$
with distinct colors, we see that every two vertices of $G$ are
connected by $k$ internally disjoint rainbow paths and that
$rc_k(G)$ is defined for every $1\leq k\leq \kappa$. So $rc_k(G)$ is
well-defined. Furthermore, $rc_k(G)\leq rc_j(G)$ for $1\leq k\leq
j\leq \kappa$. Note that this new defined $rainbow$
$k$-$connectivity$ computes the number of $colors$, this is distinct
with connectivity (edge-connectivity) which computes the number of
internally (edge) disjoint $paths$. We can also call it $rainbow$
$k$-$connection$ $number$.

Now we introduce another generalization of rainbow connection number
by Chartrand, Okamoto and Zhang \cite{Chartrand 4}. Let $G$ be an
edge-colored nontrivial connected graph of order $n$. A tree $T$ in
$G$ is a $rainbow$ $tree$ if no two edges of $T$ are colored the
same. Let $k$ be a fixed integer with $2\leq k\leq n$. An edge
coloring of $G$ is called a $k$-$rainbow~coloring$ if for every set
$S$ of $k$ vertices of $G$, there exists a rainbow tree in $G$
containing the vertices of $S$. The $k$-$rainbow~index$
${r\textsl{x}}_k(G)$ of $G$ is the minimum number of colors needed
in a $k$-rainbow coloring of $G$. A $k$-rainbow coloring using
${r\textsl{x}}_k(G)$ colors is called a $minimum$
$k$-$rainbow~coloring$. Thus ${r\textsl{x}}_2(G)$ is the rainbow
connection number $rc(G)$ of $G$. It follows, for every nontrivial
connected graph $G$ of order $n$, that ${r\textsl{x}}_2(G)\leq
{r\textsl{x}}_3(G)\leq \cdots\leq {r\textsl{x}}_n(G)$.

The above four new graph-parameters are all defined in edge-colored
graphs. Krivelevich and Yuster \cite{Krivelevich} naturally
introduced a new parameter corresponding to rainbow connection
number which is defined on vertex-colored graphs. A vertex-colored
graph $G$ is $rainbow~vertex$-$connected$ if any two vertices are
connected by a path whose $internal$ vertices have distinct colors.
A vertex-coloring under which $G$ is rainbow vertex-connected is
called a $rainbow$ $vertex$-$coloring$. The
$rainbow$~$vertex$-$connection$~$number$ of a connected graph $G$,
denoted by $rvc(G)$, is the smallest number of colors that are
needed in order to make $G$ rainbow vertex-connected. The minimum
rainbow vertex-coloring is defined similarly. Obviously, we always
have $rvc(G)\leq n-2$ (except for the singleton graph), and
$rvc(G)=0$ if and only if $G$ is a clique. Also clearly, $rvc(G)\geq
diam(G)-1$ with equality if the diameter of $G$ is 1 or 2.

Note that $rvc(G)$ may be much smaller than $rc(G)$ for some graph
$G$. For example, $rvc(K_{1,n-1})=1$ while $rc(K_{1,n-1})=n-1$.
$rvc(G)$ may also be much larger than $rc(G)$ for some graph $G$.
For example, take $n$ vertex-disjoint triangles and, by designating
a vertex from each of them, add a complete graph on the designated
vertices. This graph has $n$ cut-vertices and hence $rvc(G)\geq n$.
In fact, $rvc(G)=n$ by coloring only the cut-vertices with distinct
colors. On the other hand, it is not difficult to see that
$rc(G)\leq 4$. Just color the edges of the $K_n$ with, say, color 1,
and color the edges of each triangle with the colors $2,3,4$.

In Section 2, we will focus on the rainbow connection number and
strong rainbow connection number. We collect many upper bounds for
these two parameters. From Section 3 to Section 5, we survey on the
other three parameters: rainbow $k$-connectivity, $k$-rainbow index,
rainbow vertex-connection number, respectively. In the last section,
we sum up the results on algorithms and computational complexity.

\subsection{Terminology and notations}

All graphs considered in this survey are finite, simple and
undirected. We follow the notations and terminology of \cite{Bondy}
for all those not defined here. We use~$V(G)$ and $E(G)$ to denote
the set of vertices and the set of edges of $G$, respectively. For
any subset $X$ of $V(G)$, let $G[X]$ denote the subgraph induced by
$X$, and $E[X]$ the edge set of $G[X]$; similarly, for any subset
$F$ of $E(G)$, let $G[F]$ denote the subgraph induced by $F$. Let
$\mathcal {G}$ be a set of graphs, then $V(\mathcal
{G})=\bigcup_{G\in \mathcal {G}}{V(G)}$, $E(\mathcal
{G})=\bigcup_{G\in \mathcal {G}}{E(G)}$. We define a $clique$ in a
graph $G$ to be a complete subgraph of $G$, and a $maximal~clique$
is a clique that is not contained in any larger clique of $G$. For a
set $S$, $|S|$ denotes the cardinality of $S$. An edge in a
connected graph is called a $bridge$, if its removal disconnects the
graph. A graph with no bridges is called a $bridgeless$ graph. A
vertex is called $pendant$ if its degree is 1. We call a path of $G$
with length $k$ a $pendant$~$k$-$length$~$path$ if one of its end
vertex has degree 1 and all inner vertices have degree 2 in $G$. By
definition, a pendant $k$-length path contains a pendant
$\ell$-length path ($1\leq \ell \leq k$). A pendant 1-length path is
a $pendant$~$edge$. We denote $C_n$ a cycle with $n$ vertices. For
$n\geq 3$, the $wheel$ $W_n$ is constructed by joining a new vertex
to every vertex of $C_n$. We use $g(G)$ to denote the girth of $G$,
that is, the length of a shortest cycle of $G$.

Let $G$ be a connected graph. Recall that the distance between two
vertices $u$ and $v$ in $G$, denoted by $d(u, v)$, is the length of
a shortest path between them in $G$. The $eccentricity$ of a vertex
$v$ is $ecc(v):=\max_{x\in V(G)}{d(v,x)}$. The $diameter$ of $G$ is
$diam(G):=\max_{x\in V(G)}{ecc(x)}$. The $radius$ of $G$ is
$rad(G):=\min_{x\in V(G)}{ecc(x)}$. Distance between a vertex $v$
and a set $S \subseteq V(G)$ is $d(v,S):=\min_{x\in S}{d(v,x)}$. The
$k$-$step~open~neighbourhood$ of a set $S \subseteq V(G)$ is
$N^k(S):=\{x\in V(G)|d(x,S)=k\}$, $k\in \{0,1,2,\cdots \}$. A set $D
\subseteq V(G)$ is called a $k$-$step$ $dominating$ $set$ of $G$, if
every vertex in $G$ is at a distance at most $k$ from $D$. Further,
if $D$ induces a connected subgraph of G, it is called a $connected$
$k$-$step$ $dominating$ $set$ of $G$. The cardinality of a minimum
connected $k$-step dominating set in $G$ is called its $connected$
$k$-$step$ $domination$ $number$, denoted by $\gamma_c^k(G)$. We
call a two-step dominating set $k$-$strong$ \cite{Krivelevich} if
every vertex that is not dominated by it has at least $k$ neighbors
that are dominated by it. In \cite{Chandran}, Chandran, Das,
Rajendraprasad and Varma made two new definitions which will be
useful in the sequel. A dominating set $D$ in a graph $G$ is called
a $two$-$way~ dominating~ set$ if every pendant vertex of $G$ is
included in $D$. In addition, if $G[D]$ is connected, we call $D$ a
$connected~ two$-$way~ dominating~ set$. A (connected) two-step
dominating set $D$ of vertices in a graph $G$ is called a
$(connected)~ two$-$way~ two$-$step~ dominating~ set$ if $(i)$ every
pendant vertex of $G$ is included in $D$ and $(ii)$ every vertex in
$N^2(D)$ has at least two neighbours in $N^1(D)$. Note that if
$\delta(G)\geq 2$, then every (connected) dominating set in $G$ is a
(connected) two-way dominating set.

A subgraph $H$ of a graph $G$ is called $isometric$ if distance
between any pair of vertices in $H$ is the same as their distance in
$G$. The size of a largest isometric cycle in $G$ is denoted by
$iso(G)$. A graph is called $chordal$ if it contains no induced
cycles of length greater than 3. The $chordality$ of a graph $G$ is
the length of a largest induced cycle in $G$. Note that every
isometric cycle is induced and hence $iso(G)$ is at most the
chordality of $G$. For $k\leq \alpha(G)$, we use $\sigma_k(G)$ to
denote the minimum degree sum that is taken over all independent
sets of $k$ vertices of $G$, where $\alpha(G)$ is the number of
elements of an maximum independent set of $G$.

\section{(Strong) Rainbow connection number}
\subsection{Basic results} In \cite{Chartrand 1}, Chartrand, Johns, McKeon and Zhang
did some basic research on the (strong) rainbow connection numbers
of graphs. They determined the precise (strong) rainbow connection
numbers of several special graph classes including trees, complete
graphs, cycles, wheel graphs, complete bipartite graphs and complete
multipartite graphs.
\begin{pro}\label{thm101}\cite{Chartrand 1}
Let $G$ be a nontrivial connected graph of size $m$. Then\\
(a)~ $rc(G)=1$ if and only if $G$ is complete,~$src(G)=1$ if and only if $G$ is complete;\\
(b)~ $rc(G)=2$ if and only if $src(G)=2$;\\
(c)~ $rc(G)=m$ if and only if $G$ is a tree, $src(G)=m$ if and only
if $G$ is a tree.
\end{pro}

\begin{pro}\label{thm102}\cite{Chartrand 1}
For each integer $n\geq 4$, $rc(C_n)=src(C_n)=\lceil
\frac{n}{2}\rceil$.
\end{pro}

\begin{pro}\label{thm103}\cite{Chartrand 1}
For each integer $n\geq 3$, we have

\[
 rc(W_n)=\left\{
   \begin{array}{ll}
     1 &\mbox {if~$n=3$,}\\
     2 &\mbox {if~$4\leq n\leq 6$,}\\
     3 &\mbox {if~$n\geq 7$.}
   \end{array}
   \right.
\] and $src(W_n)=\lceil \frac{n}{3}\rceil$.

\end{pro}

\begin{pro}\label{thm104}\cite{Chartrand 1}
For integers $s$ and $t$ with $2\leq s\leq t$,
$rc(K_{s,t})=\min\{\lceil \sqrt[s]{t}\rceil,4\},$ and for integers
$s$ and $t$ with $1\leq s\leq t$, $src(K_{s,t})=\lceil
\sqrt[s]{t}\rceil.$
\end{pro}

\begin{pro}\label{thm105}\cite{Chartrand 1}
Let $G=K_{n_1,n_2,\ldots,n_k}$ be a complete $k$-partite graph,
where $k\geq 3$ and $n_1\leq n_2\leq \ldots \leq n_k$ such that
$s=\sum_{i=1}^{k-1}{n_i}$ and $t=n_k$. Then

\[
 rc(G)=\left\{
   \begin{array}{ll}
     1 &\mbox {if~$n_k=1$,}\\
     2 &\mbox {if~$n_k\geq 2$~and~$s>t$,}\\
     \min\{\lceil \sqrt[s]{t}\rceil,3\} &\mbox {if~$s\leq t$.}
   \end{array}
   \right.
\] and
\[
 src(G)=\left\{
   \begin{array}{ll}
     1 &\mbox {if~$n_k=1$,}\\
     2 &\mbox {if~$n_k\geq 2$~and~$s>t$,}\\
     \lceil \sqrt[s]{t}\rceil &\mbox {if~$s\leq t$.}
   \end{array}
   \right.
\]
\end{pro}

By Proposition \ref{thm101}, it follows that for every positive
integer $a$ and for every tree $T$ of size $a$, $rc(T)=src(T)=a$.
Furthermore, for $a\in \{1,2\}$, $rc(G)=a$ if and only if
$src(G)=a$. If $a=3,b\geq 4$, then by Proposition \ref{thm103},
$rc(W_{3b})=3$ and $src(W_{3b})=b$. For $a\geq 4$, we have the
following.

\begin{thm}\label{thm106}\cite{Chartrand 1} Let $a$ and $b$ be positive integers with $a \geq 4$ and
$b \geq \frac{5a-6}{3}$. Then there exists a connected graph $G$
such that $rc(G)=a$ and $src(G)=b$.
\end{thm}

Then, combining Propositions \ref{thm101} and \ref{thm103} with
Theorem \ref{thm106}, they got the following result.

\begin{cor}\label{thm107}\cite{Chartrand 1}
Let $a$ and $b$ be positive integers. If $a = b$ or $3 \leq a < b$
and $b \geq \frac{5a-6}{3}$ , then there exists a connected graph
$G$ such that $rc(G) = a$ and $src(G) = b$.
\end{cor}

Finally, they thought the question that whether the condition $b
\geq \frac{5a-6}{3}$ can be deleted ? and raised the following
conjecture:

\begin{con}\label{thm108}\cite{Chartrand 1}
Let $a$ and $b$ be positive integers. Then there exists a connected
graph $G$ such that $rc(G) = a$ and $src(G) = b$ if and only if
$a=b\in \{1,2\}$ or $3\leq a\leq b$.
\end{con}

In \cite{Chen-Li}, Chen and Li gave a confirmative solution to this
conjecture by showing a class of graphs with given rainbow
connection number $a$ and strong rainbow connection number $b$.

From the above several propositions, we know $rc(G)=src(G)$ hold for
some special graph classes. A difficult problem is following:

\begin{prob}\label{prob1}
Characterize graphs $G$ for which $rc(G)=src(G)$, or, give some
sufficient conditions to guarantee $rc(G)=src(G)$.
\end{prob}

Recall the fact that if $H$ is a connected spanning subgraph of a
nontrivial (connected) graph $G$, then $rc(G)\leq rc(H)$. This fact
is very useful to bounding the value of $rc(G)$ by giving bounds for
its connected spanning subgraphs. We have noted that if, in
addition, $diam(H)=2$, then $src(G)\leq src(H)$. The authors of
\cite{Chartrand 1} naturally raised the following conjecture:

\begin{con}\label{thm109}\cite{Chartrand 1} If $H$ is a connected spanning subgraph of a
nontrivial (connected) graph $G$, then $src(G)\leq src(H)$.
\end{con}

\begin{figure}[!hbpt]
\begin{center}
\includegraphics[scale=1.0000000]{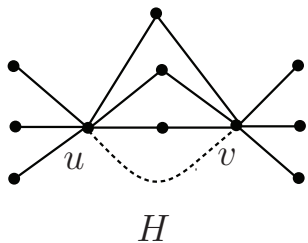}
\end{center}
\caption{A counterexample to Conjecture
\ref{thm109}.}\label{figure2}
\end{figure}

Recently, this conjecture was disproved by Chakraborty, Fischer,
Matsliah and Yuster \cite{Chakraborty}. They showed the following
example: see Figure \ref{figure2}, here $G$ is obtained from $H$ by
adding the edge $e=uv$, then $H$ is a connected spanning subgraph of
$G$. It is easy to show that there is a strong rainbow coloring of
$H$ which costs six colors, but the graph $G$ costs at least seven
colors to ensure its strong rainbow connection.

Suppose that $G$ contains two bridges $e=uv$ and $f=xy$. Then
$G-e-f$ contains three components $G_i(1\leq i\leq 3)$, where two of
these components contain one of $u,v,x$ and $y$ and the third
component contains two of these four vertices, say $u\in V(G_1)$,
$x\in V(G_2)$ and $v,y\in V(G_3)$. If $S$ is a set of $k$ vertices
contains $u$ and $x$, then every tree whose vertex set contains $S$
must also contain the edges $e$ and $f$. This gives us a necessary
condition for an edge-colored graph to be $k$-rainbow colored.
\begin{obs}\label{ob1}\cite{Chartrand 4}
Let $G$ be a connected graph of order $n$ containing two bridges $e$
and $f$. For each integer $k$ with $2\leq k\leq n$, every
$k$-rainbow coloring of $G$ must assign distinct colors to $e$ and
$f$.
\end{obs}

From Observation \ref{ob1}, we know that if $G$ is rainbow connected
under some edge-coloring, then any two bridges obtain distinct
colors.

\subsection{Upper bounds for rainbow connection number}

We know that it is almost impossible to give the precise rainbow
connection number of a given arbitrary graph, so we aim to give some
nice bounds for it, especially sharp upper bounds.

 In \cite{Caro}, Caro, Lev, Roditty, Tuza and Yuster investigated the extremal
graph-theoretic behavior of rainbow connection number. Motivated by
the fact that there are graphs with minimum degree 2 and with
$rc(G)=n-3$ (just take two vertex-disjoint triangles and connect
them by a path of length $n-5$), it is interesting to study the
rainbow connection number of graphs with minimum degree at least 3
and they thought of the following question: is it true that minimum
degree at least 3 guarantees $rc(G)\leq {\alpha}n $ where $\alpha <
1$ is independant of $n$? This turns out to be true, and they
proved:

\begin{thm}\label{thm110}\cite{Caro}
If $G$ is a connected graph with $n$ vertices and $\delta(G)\geq 3$,
then $rc(G)< \frac{5}{6}{n}$.
\end{thm}

In the proof of Theorem \ref{thm110}, they first gave an upper bound
for the rainbow connection number of 2-connected graphs (see Theorem
\ref{thm118}), then from it, they next derived an upper bound for
the rainbow connection number of connected bridgeless graphs (see
Theorem \ref{thm120}).

The constant 5/6 appearing in Theorem \ref{thm110} is not optimal,
but it probably cannot be replaced with a constant smaller than
$\frac{3}{4}$, since there are 3-regular connected graphs with
$rc(G)=diam(G)=\frac{3n-10}{4}$, and one of such graphs can be
constructed as follows \cite{Ingo}: Take two vertex disjoint copies
of the graph $K_5-P_3$ and label the two vertices of degree 2 with
$w_1$ and $w_{2k+2}$, where $k \geq 1$ is an integer. Next join
$w_1$ and $w_{2k+2}$ by a path of length $2k+1$ and label the
vertices with $w_1, w_2, \cdots , w_{2k+2}$. Now for $1\leq i\leq k$
every edge $w_{2i}w_{2i+1}$ is replaced by a $K_4-e$ and we identify
the two vertices of degree 2 in $K_4-e$ with $w_{2i}$ and
$w_{2i+1}$. The resulting graph $G_{4k+10}$ is 3-regular, has order
$n=4k+10$ and $rc(G_{4k+10})=diam(G_{4k+10})=3k+5=\frac{3n-10}{4}$.
Then Caro, Lev, Roditty, Tuza and Yuster conjectured:

\begin{con}\label{thm111}\cite{Caro} If $G$ is a connected graph with $n$ vertices and $\delta(G)\geq
3$, then $rc(G)< \frac{3}{4}{n}$.
\end{con}

Schiermeyer proved the conjecture in \cite{Ingo} by showing the
following result:

\begin{thm}\label{thm112}\cite{Ingo} If $G$ is a connected graph with $n$ vertices and $\delta(G)\geq
3$, then $rc(G)< \frac{3n-1}{4}$.
\end{thm}

For 2-connected graphs Theorem \ref{thm112} is true by Theorem
\ref{thm118}. Hence it remains to prove it for graphs with
connectivity 1. Schiermeyer extended the concept of rainbow
connection number as follows: Additionally we require that any two
edges of $G$ have different colors whenever they belong to different
blocks of $G$. The corresponding rainbow connection number will be
denoted by $rc^*(G)$. Then they derived Theorem \ref{thm112} by
first proving the following result: let $G$ be a connected graph
with $n$ vertices, connectivity 1, and $\delta\geq 3$, then
$rc^*(G)\leq \frac{3n-10}{4}$.

Not surprisingly, as the minimum degree increases, the graph would
become more dense and therefore the rainbow connection number would
decrease. Specifically, Caro, Lev, Roditty, Tuza and Yuster also
proved the following upper bounds in term of minimum degree.

\begin{thm}\label{thm113}\cite{Caro}
If $G$ is a connected graph with $n$ vertices and minimum degree
$\delta$, then $$rc(G)\leq
\min\{n\frac{\ln{\delta}}{\delta}(1+o_{\delta}(1)),
n\frac{4\ln{\delta}+3}{\delta}\}.$$
\end{thm}

In the proof, they used the concept of a connected two-dominating
set (A set of vertices $S$ of $G$ is called a $connected$
$two$-$dominating$ $set$ if $S$ induces a connected subgraph of $G$
and, furthermore, each vertex outside of $S$ has at least two
neighbours in $S$) and the probabilistic method. They showed that in
any case it cannot be improved below
$\frac{3n}{\delta+1}-\frac{\delta+7}{\delta+1}$ as they constructed
a connected $n$-vertex graph with minimum degree $\delta$ and this
diameter: Take $m$ copies of $K_{\delta+1}$, denoted by $X_1,\cdots,
X_m$ and label the vertices of $X_i$ with $x_{i,1},\cdots,
x_{i,\delta+1}$. Take two copies of $K_{\delta+2}$, denoted by $X_0,
X_{m+1}$ and similarly label their vertices. Now, connect $x_{i,2}$
with $x_{i+1,1}$ for $i=0, \cdots, m$ with an edge, and delete the
edges $(x_{i,1}, x_{i,2})$ for $i= 0, \cdots, m +1$. The obtained
graph has $n =(m + 2)(\delta+ 1)+ 2$ vertices, and minimum degree
$\delta$ (and maximum degree $\delta +1$). It is straightforward to
verify that a shortest path from $x_{0,1}$ to $x_{m+1,2}$ has length
$3m+ 5= \frac{3n}{\delta+1}-\frac{\delta+7}{\delta+1}$.

This, naturally, raised the open problem of determining the true
behavior of $rc(G)$ as a function of $\delta$.

In \cite{Chakraborty}, Chakraborty, Fischer, Matsliah and Yuster
proved that any connected $n$-vertex graph with minimum degree
$\Theta(n)$ has a bounded rainbow connection.

\begin{thm}\label{thm114}\cite{Chakraborty}
For every $\epsilon > 0$ there is a constant $C = C(\epsilon)$ such
that if $G$ is a connected graph with $n$ vertices and minimum
degree at least ${\epsilon}n$, then $rc(G) \leq C$.
\end{thm}

The proof of Theorem \ref{thm114} is based upon a modified
degree-form version of Szemer\'{e}di Regularity Lemma (see
\cite{komlos} for a good survey on Regularity Lemma) that they
proved.

The above lower bound construction suggests that the logarithmic
factor in their upper bound may not be necessary and that, in fact
$rc(G)\leq Cn/{\delta}$ where $C$ is a universal constant. If true,
notice that for graphs with a linear minimum degree ${\epsilon}n$,
this implies that $rc(G)$ is at most $C/{\epsilon}$. However,
Theorem \ref{thm114} does not even guarantee the weaker claim that
$rc(G)$ is a constant. The constant $C = C(\epsilon)$ they obtained
is a $tower~function$ in $1/{\epsilon}$ and in particular extremely
far from being reciprocal to $1/{\epsilon}$.

Finally, Krivelevich and Yuster in \cite{Krivelevich} determined the
behavior of $rc(G)$ as a function of $\delta(G)$ and resolved the
above-mentioned open problem.

\begin{thm}\label{thm115}\cite{Krivelevich}
A connected graph $G$ with $n$ vertices has
$rc(G)<\frac{20n}{\delta(G)}$.
\end{thm}

The proof of Theorem \ref{thm115} uses the concept of connected
two-step dominating set. Krivelevich and Yuster first proved that
for a connected graph $H$ with minimum degree $k$ and $n$ vertices,
there exists a two-step dominating set $S$ whose size is at most
$\frac{n}{k+1}$ and there is a connected two-step dominating set
$S'$ containing $S$ with $|S'|\leq 5|S|-4$. They found two
edge-disjoint spanning subgraphs in a graph $G$ with minimum degree
at least $\lfloor \frac{\delta-1}{2}\rfloor$. Then they derived a
rainbow coloring for $G$ by giving a rainbow coloring to each
subgraphs according to its connected two-step dominating set.

The authors noted that the constant 20 obtained by their proof is
not optimal and can be slightly improved with additional effort.
However, from the example below Theorem \ref{thm113}, one cannot
expect to replace $C$ by a constant smaller than 3.

Motivated by the results of Theorems \ref{thm112}, \ref{thm113} and
\ref{thm115},  Schiermeyer raised the following open problem in
\cite{Ingo}.

\begin{prob}\label{prob2}\cite{Ingo}
 For every $k \geq 2$ find a minimal constant
$c_k$ with $0 < c_k \leq 1$ such that $rc(G) \leq {c_k}n$ for all
graphs $G$ with minimum degree ${\delta}(G) \geq k$. Is it true that
$c_k=\frac{3}{k+1}$ for all $k \geq 2$ $?$
\end{prob}

This is true for $k = 2,3$ as shown before ($c_2 = 1$ and $c_3 =
\frac{3}{4})$.

Recently, Chandran, Das, Rajendraprasad and Varma \cite{Chandran}
nearly settled the above problem. They used the concept of a
connected two-way two-step dominating set in the argument and they
first proved the following result.

\begin{thm}\label{thm116}\cite{Chandran}
If $D$ is a connected two-way two-step dominating set in a graph
$G$, then $rc(G) \leq rc(G[D]) + 6$.
\end{thm}

Furthermore, they gave a nearly sharp bound for the size of $D$ by
showing that every connected graph $G$ of order $n \geq 4$ and
minimum degree $\delta$ has a connected two-way two-step dominating
set $D$ of size at most $\frac{3n}{\delta+1}-2$; moreover, for every
$\delta \geq 2$, there exist infinitely many connected graphs $G$
such that $\gamma_c^2(G)\geq \frac{3(n-2)}{\delta+1}-4$. Then the
following result is easy.

\begin{thm}\label{thm117}\cite{Chandran}
For every connected graph $G$ of order $n$ and minimum degree
$\delta$, $$rc(G)\leq \frac{3n}{\delta+1}+3.$$ Moreover, for every
$\delta \geq 2$, there exist infinitely many connected graphs $G$
such that $rc(G)\geq \frac{3(n-2)}{\delta+1}-1$.
\end{thm}

Theorem \ref{thm117} answers Problem \ref{prob2} in the affirmative
but up to an additive constant of 3. Moreover, this bound is seen to
be tight up to additive factors by the construction mentioned in
\cite{Caro} (see the example below Theorem \ref{thm113}) and
\cite{Erdos}. And therefore, for graphs with linear minimum degree
${\epsilon}n$, the rainbow connection number is bounded by a
constant.

Recently, Dong and Li \cite{Dong-Li} derived an upper bound on
rainbow connection numbers of graphs under given degree sum
condition $\sigma_2$. Recall that for a graph $G$, $\sigma_2(G)=\min
\{d(u) + d(v)~|~u, v ~are~independent~in~G \}$. Clearly, the degree
sum condition $\sigma_2$ is weaker than the minimum degree
condition.

\begin{thm}\label{thm155}\cite{Dong-Li} For a connected graph $G$ of order
$n$, $rc(G)\leq 6{\frac{n-2}{\sigma_2+2}} +7$.
\end{thm}

Similar to the method of Theorem \ref{thm117}, they derived that
every connected graph $G$ of order $n$ with at most one pendant
vertex has a connected two-way two-step dominating set $D$ of size
at most $6{\frac{n-2}{\sigma_2+2}} +2$. Then by using Theorem
\ref{thm116}, they got the theorem.\emph{\emph{}}

From the example below Theorem \ref{thm113}, we know their bound are
seen to be tight up to additive factors. Note that by the definition
of $\sigma_2$, we know $\sigma_2\geq 2\delta$, so from Theorem
\ref{thm155}, we can derive $rc(G)\leq 6{\frac{n-2}{\sigma_2+2}}
+7\leq \frac{3(n-2)}{\delta+1}+7$. And the bound in Theorem
\ref{thm155} can be seen as an improvement of that in Theorem
\ref{thm117}.

 With respect to the the relation between $rc(G)$ and the
connectivity $\kappa(G)$, mentioned in \cite{Ingo}, Broersma asked a
question at the IWOCA workshop:

\begin{prob}\label{prob4}\cite{Ingo}
What happens with the value $rc(G)$ for graphs with higher
connectivity?
\end{prob}

For $\kappa(G)=1$, Theorem \ref{thm112} means that if $G$ is a graph
of order $n$, connectivity $\kappa(G)=1$ and $\delta \geq 3$. Then
$rc(G)\leq \frac{3n-1}{4}$. For $\kappa(G)=2$, in the proof of
Theorem \ref{thm110}, as we mentioned above, Caro, Lev, Roditty,
Tuza and Yuster derived:

\begin{thm}\label{thm118}\cite{Caro}
If $G$ is a $2$-connected graph with $n$ vertices then $rc(G)\leq
\frac{2n}{3}$.
\end{thm}

That is, if $G$ is a graph of order $n$, connectivity $\kappa(G)=2$.
Then $rc(G)\leq \frac{2n}{3}$.

From Theorem \ref{thm117}, we can easily obtain an upper bound of
the rainbow connection number according to the connectivity:
$$rc(G)\leq \frac{3n}{\delta+1}+3\leq \frac{3n}{\kappa+1}+3.$$
Therefore, for $\kappa(G)=3$, $rc(G)\leq \frac{3n}{4}+3$; for
$\kappa(G)=4$, $rc(G)\leq \frac{3n}{5}+3$. Motivated by the results
in \cite{Caro}, and by using the Fan Lemma, Li and Shi
\cite{Li-Shi2} improved this bound by showing the following result.

\begin{thm}\label{thm119}(\cite{Li-Shi2})
If $G$ is a 3-connected graph with $n$ vertices, then $rc(G) \leq
\frac{3(n+1)}{5}$.
\end{thm}

However, for general connectivity, there is no upper bound which is
better than $\frac{3n}{\kappa+1}+3$.

The following result is an important ingredient in the proof of
Theorem \ref{thm110} in \cite{Caro}.

\begin{thm}\label{thm120}\cite{Caro}
If $G$ is a connected bridgeless graph with $n$ vertices, then
$rc(G) \leq \frac{4n}{5}-1$.
\end{thm}

From Theorem \ref{thm117}, we can also easily obtain an upper bound
of the rainbow connection number according to the edge-connectivity
$\lambda$:
$$rc(G)\leq \frac{3n}{\delta+1}+3\leq \frac{3n}{\lambda+1}+3.$$
Note that all the above upper bounds are determined by $n$ and other
parameters such as (edge)-connectivity, minimum degree. Diameter of
a graph, and hence its radius, are obvious lower bounds for rainbow
connection number. Hence it is interesting to see if there is an
upper bound which is a function of the radius $r$ or diameter alone.
Such upper bounds were shown for some special graph classes in
\cite{Chandran} which we will introduce in the sequel. But, for a
general graph, the rainbow connection number cannot be upper bounded
by a function of r alone. For instance, the star $K_{1,n}$ has a
radius 1 but rainbow connection number $n$. Still, the question of
whether such an upper bound exists for graphs with higher
connectivity remains. Basavaraju, Chandran, Rajendraprasad and
Ramaswamy \cite{Basavaraju} answered this question in the
affirmative. The key of their argument is the following lemma, and
in the proof of this lemma, we can obtain a connected $(k-1)$-step
dominating set from a connected $k$-step dominating set.
\begin{lem}\label{thm153}\cite{Basavaraju}
If $G$ is a bridgeless graph, then for every connected $k$-step
dominating set $D^k$ of $G$, $k\geq 1$, there exists a connected
$(k-1)$-step dominating set $D^{k-1}\supset D^k$ such that
$$rc(G[D^{k-1}])\leq rc(G[D^k])+\min\{2k+1, \zeta\},$$ where
$\zeta=iso(G)$.
\end{lem}

Given a graph $G$ and a set $D\subset V(G)$, A $D$-$ear$ is a path
$P=(x_0,x_1,\cdots,x_p)$ in $G$ such that $P \cap D=\{x_0,x_p \}$.
$P$ may be a closed path, in which case $x_0=x_p$. Further, $P$ is
called an $acceptable$ $D$-ear if either $P$ is a shortest $D$-ear
containing $(x_0,x_1)$ or $P$ is a shortest $D$-ear containing
$(x_{p-1},x_p)$. Let $\mathcal {A}=\{a_1,a_2,\cdots\}$ and $\mathcal
{B}=\{b_1,b_2,\cdots\}$ be two pools of colors, none of which are
used to color $G[D^k]$. A $D^k$-$ear$ $P=(x_0,x_1,\cdots,x_p)$ will
be called $evenly~colored$ if its edges are colored
$a_1,a_2,\cdots,a_{\lceil \frac{p}{2}\rceil},b_{\lfloor
\frac{p}{2}\rfloor},\cdots,b_2,b_1$ in that order. Basavaraju,
Chandran, Rajendraprasad and Ramaswamy proved this lemma by
constructing a sequence of sets  $D^k=D_0\subset D_1\subset \cdots
\subset D_t=D^{k-1}$ and coloring the new edges in every induced
graph $G[D_i]$ such that the following property is maintained for
all $0 \leq i\leq t$: every $x\in D_i\backslash D^k$ lies in an
evenly colored acceptable $D^k$-ear in $G[D_i]$.

The following theorem can be derived from Lemma \ref{thm153} easily.
\begin{thm}\label{thm121}\cite{Basavaraju}
For every connected bridgeless graph $G$,
$$rc(G)\leq \sum_{i=1}^r{\min}\{2i+1,\zeta\}\leq r{\zeta},$$
where $r$ is the radius of $G$.
\end{thm}

Theorem \ref{thm121} has two corollaries.

\begin{cor}\label{thm122}\cite{Basavaraju}
For every connected bridgeless graph $G$ with radius $r$,
$$rc(G)\leq r(r+2).$$
Moreover, for every integer $r\geq 1$, there exists a bridgeless
graph with radius $r$ and $rc(G)=r(r+2)$.
\end{cor}

\begin{cor}\label{thm123}\cite{Basavaraju}
For every connected bridgeless graph $G$ with radius $r$ and
chordality $k$,
$$rc(G)\leq \sum_{i=1}^r{\min}\{2i+1,k\}\leq rk.$$
Moreover, for every two integers $r\geq 1$ and $3\leq k\leq 2r+1$,
there exists a bridgeless graph $G$ with radius $r$ and chordality
$k$ such that $rc(G)=\sum_{i=1}^r{\min}\{2i+1,k\}$.
\end{cor}

Corollary \ref{thm122} answered the above question in the
affirmative, the bound is sharp and is a function of the radius $r$
alone. Basavaraju, Chandran, Rajendraprasad and Ramaswamy also
demonstrated that the bound cannot be improved even if we assume
stronger connectivity by constructing a $\kappa$-vertex-connected
graph of radius $r$ whose rainbow connection number is $r(r+2)$ for
any two given integers $\kappa, r \geq 1$: Let $s(0):= 0,
s(i):=2{\sum_{j=r}^{r-i+1}{j}}$ for $1\leq i\leq r$ and $t:=
s(r)=r(r+1)$. Let $V=V_0\uplus V_1\uplus \cdots \uplus V_t$ where
$V_i=\{x_{i,0}, x_{i,1},\cdots, x_{i,\kappa-1}\}$ for $0\leq i\leq
t-1$ and $V_t=\{x_{t,0}\}$. Construct a graph $X_{r,\kappa}$ on $V$
by adding the following edges. $E(X)=\{\{x_{i,j}, x_{i',j'}\}:
|i-i'|\leq 1\} \cup \{\{x_{s(i),0}, x_{s(i+1),0} \}: 0\leq i\leq
r-1\}$.

Corollary \ref{thm123} generalises a result from \cite{Chandran}
that the rainbow connection number of any bridgeless chordal graph
is at most three times its radius as the chordality of a chordal
graph is three.

In \cite{Caro}, Caro, Lev, Roditty, Tuza and Yuster also derived a
result which gives an upper bound for rainbow connection number
according to the order and the number of vertex-disjoint cycles.
Here $\chi'(G)$ is the chromatic index of $G$.
\begin{thm}\label{thm150}\cite{Caro}
Suppose $G$ is a connected graph with $n$ vertices, and assume that
there is a set of vertex-disjoint cycles that cover all but $s$
vertices of $G$. Then $rc(G)< 3n/4+s/4-1/2$. In particular:\\
$(i).$ If $G$ has a 2-factor then $rc(G)< 3n/4$.\\
$(ii).$ If $G$ is $k$-regular and $k$ is even then $rc(G)< 3n/4$.\\
$(iii).$ If $G$ is $k$-regular and $\chi'(G)=k$ then $rc(G)< 3n/4$.
\end{thm}

Another approach for achieving upper bounds is based on the size
(number of edges) $m$ of the graph. Those type of sufficient
conditions are known as Erd\H{o}s-Gallai type results. Research on
the following Erd\H{o}s-Gallai type problem has been started in
\cite{Kemnitz}.

\begin{prob}\label{prob5}\cite{Kemnitz}
For every $k$, $1\leq k\leq n-1$, compute and minimize the function
$f(n,k)$ with the following property: If $|E(G)|\geq f(n, k)$, then
$rc(G)\leq k$.
\end{prob}

In \cite{Kemnitz}, Kemnitz and Schiermeyer gave a lower bound for
$f(n,k)$, i.e.,  $f(n,k)\geq {{n-k+1}\choose{2}}+(k-1)$. They also
computed $f(n,k)$ for $k\in \{1,n-2,n-1\}$, i.e., $f(n,1)={n \choose
2}, f(n,n-1)=n-1, f(n,n-2)=n$, and obtained $f(n,2)={n-1 \choose
2}+1$ for $k=2$.

In \cite{Li-Sun7}, Li and Sun provided a new approach to investigate
the rainbow connection number of a graph $G$ according to some
constraints to its complement graph $\overline{G}$. They gave two
sufficient conditions to guarantee that $rc(G)$ is bounded by a
constant. By using the fact that $rc(G)\leq rc(H)$ where $H$ is a
connected spanning subgraph of a connected graph $G$, and the
structure of its complement graph as well as Propositions
\ref{thm104} and \ref{thm105}, they derived the following result.

\begin{thm}\label{thm124}\cite{Li-Sun7} For a connected
graph $G$, if $\overline{G}$ does not belong to the following two
cases: $(i)$~$diam(\overline{G})=2,3$, $(ii)~\overline{G}$ contains
exactly two connected components and one of them is trivial, then
$rc(G)\leq 4$. Furthermore, this bound is best possible.
\end{thm}

For the remaining cases, $rc(G)$ can be very large as discussed in
\cite{Li-Sun7}. So They add a constraint: let $\overline{G}$ be
triangle-free, then $G$ is claw-free. And they derived the following
result. In their argument, Theorem \ref{thm131} is useful.
\begin{thm}\label{thm125}\cite{Li-Sun7} For a connected graph $G$,
if $\overline{G}$ is triangle-free, then $rc(G)\leq 6$.
\end{thm}

The readers may consider the rainbow connection number of a graph
$G$ according to some other condition to its complement graph.

Chen, Li and Lian \cite{Chen-Li-Lian} investigated
Nordhaus-Gaddum-type result. A Nordhaus-Gaddum-type result is a
(sharp) lower or upper bound on the sum or product of the values of
a parameter for a graph and its complement. The name
``Nordhaus-Gaddum-type" is so given because it is Nordhaus and
Gaddum \cite{Nordhaus} who first established the following type of
inequalities for chromatic number of graphs in 1956.

\begin{thm}\label{thm126}\cite{Chen-Li-Lian} Let $G$ and $\overline{G}$ be connected with $n\geq 4$,
then $$4 \leq rc(G)+rc(\overline{G})\leq n+2.$$ Furthermore, the
upper bound is sharp for $n\geq 4$ and the low bound is sharp for
$n\geq 8$.
\end{thm}

They also proved that $rc(G)+rc(\overline{G})\geq 6$ for $n=4,5$;
and $rc(G)+rc(\overline{G})\geq 5$ for $n=6,7$ and these bounds are
best possible.

\subsection{For some graph classes}

Some graph classes, such as line graphs, have many special
properties, and by these properties we can get some interesting
results on their rainbow connection numbers in terms of some graph
parameters. For example, in \cite{Caro} Caro, Lev, Roditty, Tuza and
Yuster derived Theorem \ref{thm118} according to the
ear-decomposition of a 2-connected graph. In this subsection, we
will introduce some results on rainbow connection numbers of line
graphs, etc.

In \cite{Li-Sun1} and \cite{Li-Sun2}, Li and Sun studied the rainbow
connection numbers of line graphs in the light of particular
properties of line graphs shown in \cite{Hedetniemi} and
\cite{Hemminger}. They gave two sharp upper bounds for rainbow
connection number of a line graph and one sharp upper bound for
rainbow connection number of an iterated line graph.

Recall the $line~graph$ of a graph $G$ is the graph $L(G)$ (or
$L^1(G)$) whose vertex set $V(L(G))=E(G)$, and two vertices $e_1$,
$e_2$ of $L(G)$ are adjacent if and only if they are adjacent in
$G$. The $iterated~line~graph$ of a graph $G$, denoted by $L^2(G)$,
is the line graph of the graph $L(G)$. More generally, the
$k$-$iterated~line~graph$ $L^k(G)$ is the line graph of $L^{k-1}(G)
\ (k\geq 2)$. We also need the following new terminology.

\begin{figure}[!hbpt]
\begin{center}
\includegraphics[scale=1.0000000]{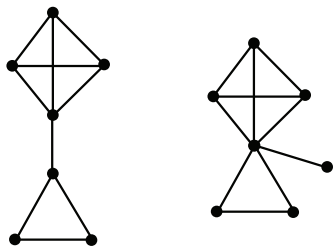}
\end{center}
\caption{A clique-forest-structure with size 6 and 2
components.}\label{figure1}
\end{figure}

For a connected graph $G$, we call $G$ a
$clique$-$tree$-$structure$, if it satisfies the following
condition: each block is a maximal clique. We call a graph $H$ a
$clique$-$forest$-$structure$, if $H$ is a disjoint union of some
clique-tree-structures, that is, each component of a
clique-forest-structure is a clique-tree-structure. By the above
condition, we know that any two maximal cliques of $G$ have at most
one common vertex. Furthermore, $G$ is formed by its maximal
cliques. The $size$ of a clique-tree(forest)-structure is the number
of its maximal cliques. An example of clique-forest-structure is
shown in Figure \ref{figure1}.
If each block of a clique-tree-structure is a triangle, we call it a
$triangle$-$tree$-$structure$. Let $\ell$ be the size of a
triangle-tree-structure. Then, by definition, it is easy to show
that there are $2\ell+1$ vertices in it. Similarly, we can give the
definition of a $triangle$-$forest$-$structure$ and there are $2\ell
+c$ vertices in a triangle-forest-structure with size $\ell$ and $c$
components. We denote $n_2$ the number of inner vertices (degrees at
least 2) of a graph.

\begin{thm}\label{thm127}\cite{Li-Sun2}
For any set $\mathcal {T}$ of $t$ edge-disjoint triangles of a
connected graph $G$, if the subgraph induced by the edge set
$E(\mathcal {T})$ is a triangle-forest-structure, then
$$rc(L(G))\leq n_2-t.$$ Moreover, the bound is sharp.
\end{thm}

\begin{thm}\label{thm128}\cite{Li-Sun2} If $G$ is a connected graph, $\mathcal {T}$
is a set of $t$ edge-disjoint triangles that cover all but $n_2'$
inner vertices of $G$ and $c$ is the number of components of the
subgraph $G[E(\mathcal {T})]$, then
$$rc(L(G))\leq t+n_2'+c.$$ Moreover, the bound is sharp.
\end{thm}

\begin{thm}\label{thm129}\cite{Li-Sun2} Let $G$ be a connected graph with $m$ edges and $m_1$
pendant 2-length paths. Then $$rc(L^2(G))\leq m-m_1.$$ The equality
holds if and only if $G$ is a path of length at least 3.
\end{thm}

In the proofs of the above three theorems, Li and Sun used the
particular structure of line graphs and the observation: If $G$ is a
connected graph and $\{E_i\}_{i\in [t]}$ is a partition of the edge
set of $G$ into connected subgraphs $G_i=G[E_i]$, then $rc(G)\leq
\sum_{i=1}^{t}{rc(G_i)}$ (see \cite{Li-Sun1}).

The above three theorems give upper bounds for rainbow connection
number of $L^k(G)(k=1,2)$ according to some parameters of $G$. One
may consider the relation between $rc(G)$ and $rc(L(G))$.

\begin{prob}\label{prob11} Determine the relationship between $rc(G)$ and
$rc(L(G))$, is there an upper bound for one of these parameters in
terms of the other?
\end{prob}

One also can consider the rainbow connection number of the general
iterated line graph $L^k(G)$ when $k$ is sufficiently large.

\begin{prob}\label{prob12} Consider the value of $rc(L^k(G))$ as $k\rightarrow \infty$, is it bounded by a constant? or, does it convergence to a function
of some graph parameters, such as the order $n$ of $G$?\qed
\end{prob}

For Problem \ref{prob12}, we know if $G$ is a cycle $C_n \ (n\geq
4)$, then $L^k(G)=G$, so $rc(L^k(G))=\lceil \frac{n}{2}\rceil$. But
for many graphs, we know, as $k$ grows, $L^k(G)$ will become more
dense, and $rc(L^k(G))$ may decrease.

An $intersection~ graph$ of a family of sets $\mathcal {F}$, is a
graph whose vertices can be mapped to the sets in $\mathcal {F}$
such that there is an edge between two vertices in the graph if and
only if the corresponding two sets in $\mathcal {F}$ have a
non-empty intersection. An $interval~ graph$ is an intersection
graph of intervals on the real line. A $unit~ interval~ graph$ is an
intersection graph of unit length intervals on the real line. A
$circular~ arc~ graph$ is an intersection graph of arcs on a circle.
An independant triple of vertices $x, y, z$ in a graph $G$ is an
$asteroidal~ triple~ (AT)$, if between every pair of vertices in the
triple, there is a path that does not contain any neighbour of the
third. A graph without asteroidal triples is called an $AT$-$free~
graph$ \cite{Corneil}. A graph $G$ is a $threshold~ graph$, if there
exists a weight function $w : V (G) \rightarrow \mathbb{R}$ and a
real constant $t$ such that two vertices $u, v \in V(G)$ are
adjacent if and only if $w(u) + w(v) \geq t$. A bipartite graph
$G(A,B)$ is called a $chain~ graph$ if the vertices of $A$ can be
ordered as $A=(a_1, a_2, \cdots,a_k)$ such that $N(a_1)\subseteq
N(a_2)\subseteq \cdots \subseteq N(a_k)$ \cite{Yannakakis}. In
\cite{Chandran}, Chandran, Das, Rajendraprasad and Varma
investigated the rainbow connection numbers of these special graph
classes. They first showed a result concerning the connected two-way
dominating sets.

\begin{thm}\label{thm131}\cite{Chandran} If $D$ is a connected
two-way dominating set in a graph $G$, then
$$rc(G)\leq rc(G[D])+3.$$
\end{thm}

They also proved that every connected graph $G$ of order $n$ and
minimum degree $\delta$ has a connected two-step dominating set $D$
of size at most $\frac{3(n-|N^2(D)|)}{\delta+1}-2$.

From Theorem \ref{thm131}, the following result can be derived.

\begin{thm}\label{thm132}\cite{Chandran} Let $G$ be a connected graph with $\delta(G)\geq 2$.
Then\\
$(i)$ if $G$ is an interval graph, $diam(G)\leq rc(G)\leq
diam(G)+1$, in particular, if $G$ is a unit interval graph, then $rc(G)=diam(G)$;\\
$(ii)$ if $G$ is $AT$-free, $diam(G)\leq rc(G)\leq diam(G)+3$;\\
$(iii)$ if $G$ is a threshold graph, $diam(G)\leq rc(G)\leq 3$;\\
$(iv)$ if $G$ is a chain graph, $diam(G)\leq rc(G)\leq 4$;\\
$(v)$ if $G$ is a circular arc graph, $diam(G)\leq rc(G)\leq
diam(G)+4$.\\
Moreover, there exist threshold graphs and chain graphs with minimum
degree at least 2 and rainbow connection number equal to the
corresponding upper bound above. There exists an $AT$-free graph $G$
with minimum degree at least 2 and $rc(G)=diam(G)+2$, which is 1
less than the upper bound above.
\end{thm}

Recall that the concept of rainbow connection number is of great use
in transferring information of high security in multicomputer
networks. Cayley graphs are very good models that have been used in
communication networks. So, it is of significance to study the
rainbow connection numbers of Cayley graphs. Li, Li and Liu
\cite{Li-Li} investigated the rainbow connection numbers of cayley
graphs on Abelian groups.

Let $\Gamma$ be a group, and let $a\in \Gamma$ be an element. We use
$\langle a\rangle$ to denote the cyclic subgroup of $\Gamma$
generated by $a$. The number of elements of $\langle a\rangle$ is
called the $order$ of $a$, denoted by $|a|$. A pair of elements $a$
and $b$ in a group commutes if $ab=ba$. A group is $Abelian$ if
every pair of its elements commutes. A $Cayley\ graph$ of $\Gamma$
with respect to $S$ is the graph $C(\Gamma,S)$ with vertex set
$\Gamma$ in which two vertices $x$ and $y$ are adjacent if and only
if $xy^{-1}\in S$ (or equivalently, $yx^{-1}\in S$), where
$S\subseteq\Gamma\setminus{\{1\}}$ is closed under taking inverse
\cite{rotman}.

\begin{thm}\label{thm133}\cite{Li-Li}
Given an Abelian group $\Gamma$ and an inverse closed set
$S\subseteq \Gamma\setminus \{1\}$, we have the following results:

{\em  $(i)$} $rc(C(\Gamma,S))\leq min\{\sum_{a\in S^*} \lceil | a
|/2\rceil\ |\ S^*\subseteq S\ is\ a\ \ minimal\ generating\ set\ of\
\Gamma\}.$

{\em  $(ii)$} If $S$ is an inverse closed minimal generating set of
$\Gamma$, then\\
\[\sum_{a\in S^*} \lfloor | a |/2\rfloor\leq rc(C(\Gamma,S))\leq
src(C(\Gamma,S))\leq\sum_{a\in S^*} \lceil | a |/2\rceil,\]

\noindent where $S^{*}\subseteq S$ is a minimal generating set of
$\Gamma$.

Moreover, if every element $a\in S$ has an even order, then

\[
rc(C(\Gamma,S))= src(C(\Gamma,S))=\sum_{a\in S^*} |a|/2.\]
\end{thm}

They also investigated the rainbow connection numbers of recursive
circulants (see \cite{park} for an introduction to recursive
circulants).

Let $G$ be an $r$-regular graph with $n$ vertices. $G$ is said to be
$strongly~regular$ and denoted by $SRG(v, k, \lambda, \mu)$ if there
are also integers $\lambda$ and $\mu$ such that every two adjacent
vertices have $\lambda$ common neighbours and every two nonadjacent
vertices have $\mu$ common neighbours. Clearly, a strongly regular
graph with parameters $(v, k, \lambda, \mu)$ is connected if and
only if $\mu \geq 1$. In \cite{Ahadi}, Ahadi and Dehghan derived the
following result: For every connected strongly regular graph $G$,
$rc(G)\leq 600$. As each strongly regular graph is a graph with
$diam(G)=2$, from our next subsection (Theorem \ref{thm154}), we
know $rc(G)\leq 5$ if $G$ is a strongly regular graph with
parameters $(v, k, \lambda, \mu)$, other than a star \cite{Li-Li2}.
But 5 may not be the optimal upper bound, so one may consider the
following question.

\begin{que}\label{que1}\cite{Ahadi}
Determine $\max \{rc(G)| ~G~is~an~SRG\}$.
\end{que}

There are other results on some special graph classes. In
\cite{Chartrand 3}, Chartrand, Johns, McKeon and Zhang investigated
the rainbow connection numbers of cages, and in \cite {Johns},
Johns, Okamoto and Zhang investigated the rainbow connection numbers
of small cubic graphs. The details are omitted.

\subsection{For dense and sparse graphs}
For any given graph $G$, we know $1\leq rc(G)\leq src(G)\leq m$.
Here a graph $G$ is called a $dense~ graph$  if its (strong) rainbow
connection number is small, especially it is close to 1; while $G$
is called a $sparse~ graph$ if its (strong) rainbow connection
number is large, especially it is close to $m$. By Proposition
\ref{thm101}, the cases that $rc(G)=1,src(G)=1$ and $rc(G)=m,
src(G)=m$ are clear. So we want to investigate other cases.

In \cite{Caro}, Caro, Lev, Roditty, Tuza and Yuster investigated the
graphs with small rainbow connection numbers, and they gave a
sufficient condition that guarantees $rc(G)=2$.

\begin{thm}\label{thm134}\cite{Caro} Any non-complete graph with $\delta(G)\geq
n/2+\log{n}$ has $rc(G)=2$.
\end{thm}

We know that having diameter 2 is a necessary requirement for having
$rc(G)=2$, although certainly not sufficient (e.g., consider a
star). Clearly, if $\delta(G)\geq n/2$ then $diam(G)= 2$, but we do
not know if this guarantees $rc(G)= 2$. The above theorem shows that
by slightly increasing the minimum degree assumption, $rc(G)= 2$
follows.

Kemnitz and Schiermeyer \cite{Kemnitz} gave a sufficient condition
to guarantee $rc(G)=2$ according to the number of edges $m$.

\begin{thm}\label{thm151}\cite{Kemnitz} Let $G$ be a connected graph
of order $n$ and size $m$. If ${{n-1}\choose 2}+1\leq m\leq
{n\choose 2} -1$, then $rc(G)=2$.
\end{thm}

Let $G=G(n,p)$ denote, as usual \cite{Alon}, the random graph with
$n$ vertices and edge probability $p$. For a graph property $A$ and
for a function $p=p(n)$, we say that $G(n,p)$ satisfies $A$ $almost~
surely$ if the probability that $G(n,p(n))$ satisfies $A$ tends to 1
as $n$ tends to infinity. We say that a function $f(n)$ is a $sharp$
$threshold~ function$ for the property $A$ if there are two positive
constants $c$ and $C$ so that $G(n,cf(n))$ almost surely does not
satisfy $A$ and $G(n,p)$ satisfies $A$ almost surely for all $p \geq
Cf(n)$. It is well known that all monotone graph properties have a
sharp threshold function (see \cite{Bollobas} and \cite{Friedgut}).
Since having $rc(G)\leq 2$ is a monotone graph property (adding
edges does not destroy this property), it has a sharp threshold
function. The following theorem establishes it.

\begin{thm}\label{thm135}\cite{Caro} $p=\sqrt{\log{n}/n}$ is a
sharp threshold function for the graph property $rc(G(n,p))\leq 2.$
\end{thm}

Theorem \ref{thm134} asserts that minimum degree $n/2+ \log n$
guarantees $rc(G)=2$. Clearly, minimum degree $n/2-1$ does not, as
there are connected graphs with minimum degree $n/2-1$ and diameter
3 (just take two vertex-disjoint cliques of order $n/2$ each and
connect them by a single edge. It is therefore interesting to raise:

\begin{prob}\label{prob7}\cite{Caro}
Determine the minimum degree threshold that guarantees $rc(G)=2$.
\end{prob}

By Proposition \ref{thm101}, we know that the problem of considering
graphs with $rc(G)=2$ is equivalent to that of considering graphs
with $src(G)=2$.

A bipartite graph which is not complete has diameter at least 3. A
proof similar to that of Theorem \ref{thm135} gives the following
result.

\begin{thm}\label{thm136}\cite{Caro} Let $c=1/{\log (9/7)}$. If $G$
is a non-complete bipartite graph with $n$ vertices
and any two vertices in the same vertex class have at least $2c{\log
n}$ common neighbors in the other vertex class, then $rc(G)=3$.
\end{thm}

The following theorem asserts however that having diameter 2 and
only logarithmic minimum degree suffices to guarantee rainbow
connection 3.

\begin{thm}\label{thm137}\cite{Chakraborty} If $G$ is an $n$-vertex
graph with diameter 2 and minimum degree at least $8\log n$, then
$rc(G)\leq 3$.
\end{thm}

Since a graph with minimum degree $n/2$ is connected and has
diameter 2, we have as an immediate result \cite{Chakraborty}: If
$G$ is an $n$-vertex graph with minimum degree at least $n/2$ then
$rc(G)\leq 3$. We know that any graph $G$ with $rc(G)=2$ must have
$diam(G)=2$, so graphs with $rc(G)=2$ belong to the graph class with
$diam(G)=2$. By Corollary \ref{thm122}, we know that for a
bridgeless graph with $diam(G)=2$, $rc(G)\leq r(r+2)\leq 8$. As
$rc(G)$ is at least the number of bridges of $G$, so it may be very
large if the number of bridges of $G$ is sufficiently large by
Observation \ref{ob1}. So there is an interesting problem:

\begin{prob}\label{prob8}
For any bridgeless graph $G$ with $diam(G)=2$, determine the
smallest constant $c$ such that $rc(G)\leq c$.
\end{prob}

Recently, Li, Li and Liu \cite{Li-Li2} derived that $c\leq 5$ by
showing the following result:

\begin{thm}\label{thm154}\cite{Li-Li2} $rc(G)\leq 5$ if $G$ is a bridgeless graph
with diameter 2; and that $rc(G)\leq k+ 2$ if $G$ is a connected
graph with diameter 2 and $k$ bridges, where $k\geq 1$.
\end{thm}

In the proof, Li, Li and Liu derived that if $G$ is a bridgeless
graph with order $n$ and diameter 2, then it is either 2-connected,
or it has only one cut vertex $v$, furthermore, $v$ is the center of
$G$ with radius 1. They showed that 5 is almost best possible as
there are infinity many bridgeless graphs with diameter 2 whose
rainbow connection numbers are 4, however they have not found
examples of such graphs with $rc(G)=5$. The bound $k+2$ is sharp as
there are infinity graphs with diameter 2 and $k$ bridges whose
rainbow connection numbers attain this bound \cite{Li-Li2}.

\begin{figure}[!hbpt]
\begin{center}
\includegraphics[scale=1.000000]{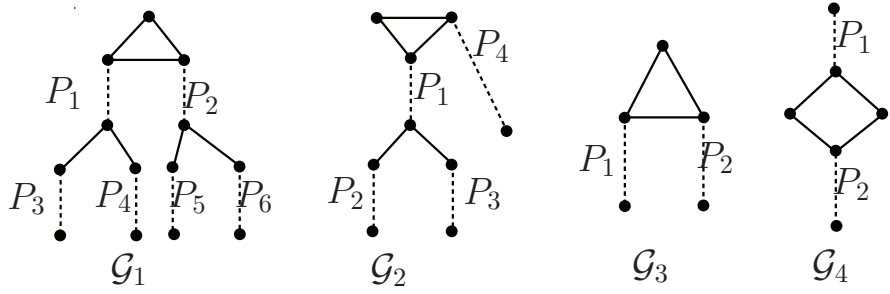}
\end{center}
\caption{The figure for the four graph classes.}\label{figure3}
\end{figure}

In \cite{Li-Sun5} and \cite{Li-Sun6}, Li and Sun investigated the
graphs with large rainbow connection numbers and strong rainbow
connection numbers, respectively. They derived the two following
results. Note that each path $P_j$ in the member of graph class
$\mathcal {G}_i \ (1\leq i\leq 4)$ of Figure \ref{figure3} may be
trivial.

\begin{thm}\label{thm138}\cite{Li-Sun5} For a connected graph $G$ with $m$ edges, we
have $rc(G)\neq m-1$; and $rc(G)=m-2$ if and only if $G$ is a
5-cycle or belongs to one of four graph classes $\mathcal {G}_i$'s
$(1\leq i\leq 4)$ shown in Figure \ref{figure3}.
\end{thm}

We now introduce two graph classes. Let $C$ be the cycle of a
unicyclic graph $G$, $V(C)=\{v_1,\cdots,v_k\}$ and $\mathcal
{T}_{G}=\{T_i:1\leq i\leq k \}$ where $T_i$ is the unique tree
containing vertex $v_i$ in subgraph $G\backslash E(C)$. We say $T_i$
and $T_j$ are $adjacent \ (nonadjacent)$ if $v_i$ and $v_j$ are
adjacent (nonadjacent) in cycle $C$. Then
let\\
$\mathcal {G}_5=\{G: G$ is a unicyclic graph, $k=3$, $\mathcal
{T}_{G}$ contains at most two nontrivial elements$\}$,\\
$\mathcal {G}_6=\{G: G$ is a unicyclic graph, $k=4$, $\mathcal
{T}_{G}$ contains two nonadjacent trivial elements and the other two
(nonadjacent) elements are paths$\}$.
\begin{thm}\label{thm139}\cite{Li-Sun6} For a connected graph $G$ with $m$ edges,
$src(G)\neq m-1$; $src(G)=m-2$ if and only if $G$ is a 5-cycle or
belongs to one of the $\mathcal {G}_i$'s $(i=5,6)$.
\end{thm}

By Proposition \ref{thm101}, Theorems \ref{thm138} and \ref{thm139},
we investigated the graphs with $rc(G)\geq m-2 \ (src(G)\geq m-2)$.
Furthermore, we have the following interesting problem.
\begin{prob}\label{prob13}
Give a sufficient condition to guarantee $rc(G)\geq
{\alpha}m~(src(G)\geq {\alpha}m)$, where $0< \alpha < 1$.
\end{prob}

\subsection{ For graph operations}
Products of graphs occur naturally in discrete mathematics as tools
in combinatorial constructions, they give rise to important classes
of graphs and deep structural problems. The extensive literature on
products that has evolved over the years presents a wealth of
profound and beautiful results \cite{Imrich}. In \cite{Li-Sun5}, Li
and Sun obtained some results on the rainbow connection numbers of
products of graphs, including Cartesian product, composition
(lexicographic product), union of graphs, etc. Actually, we know
that the line graph of a graph is also a graph operation.

We first introduce the rainbow connection number of the Cartesian
product of some graphs. The $Cartesian~product$ of graphs $G$ and
$H$ is the graph $G\square H$ whose vertex set is $V(G)\times V(H)$
and whose edge set is the set of all pairs $(u_1,v_1)(u_2,v_2)$ such
that either $u_1u_2\in E(G)$ and $v_1=v_2$, or $v_1v_2\in E(H)$ and
$u_1=u_2$. The $strong~product$ of $G$ and $H$ is the graph
$G\boxtimes H$ whose vertex set is $V(G)\times V(H)$ and whose edge
set is the set of all pairs $(u_1,v_1)(u_2,v_2)$ such that either
$u_1u_2\in E(G)$ and $v_1=v_2$, or $v_1v_2\in E(H)$ and $u_1=u_2$,
or $u_1u_2\in E(G)$ and $v_1v_2\in E(H)$. By definition, the graph
$G\square H$ is the spanning subgraph of the graph $G\boxtimes H$.
By using the definition and structure of Cartesian product, Li and
Sun derived the following.

\begin{thm}\label{thm140}\cite{Li-Sun5}
Let $G^*=G_1\square G_2\square \cdots \square G_k \ (k\geq 2)$,
where each $G_i \ (1\leq i\leq k)$ is connected. Then we have
$$rc(G^*)\leq \sum_{i=1}^k{rc(G_i)}.$$
Moreover, if $diam(G_i)=rc(G_i)$ for each $G_i$, then the equality
holds.
\end{thm}

We know that there are infinity many graph with $diam(G)=rc(G)$,
such as the unit interval graphs shown in Theorem \ref{thm132}. The
following problem could be interesting but maybe difficult:

\begin{prob}\label{prob9}
Characterize the graphs $G$ with $rc(G)=diam(G)$, or give some
sufficient conditions to guarantee that $rc(G)=diam(G)$.
\end{prob}
Similar problem for $src(G)$ can be considered.

Let $\overline{G^*}=G_1\boxtimes G_2\boxtimes \cdots \boxtimes G_k \
(k\geq 2)$, where each $G_i \ (1\leq i\leq k)$ is connected. Since
the Cartesian product of any two graphs is a spanning subgraph of
their strong product, $G^*$ is the spanning subgraph of
$\overline{G^*}$, then we have the following result.
\begin{cor}\label{thm141}\cite{Li-Sun5}
Let $\overline{G^*}=G_1\boxtimes G_2\boxtimes \cdots \boxtimes G_k \
(k\geq 2)$, where each $G_i \ (1\leq i\leq k)$ is connected. Then we
have
$$rc(\overline{G^*})\leq \sum_{i=1}^k{rc(G_i)}.$$
\end{cor}

For $i=1,2,\cdots,r$, let $m_i\geq 2$ be given integers. Consider
the graph $G$ whose vertices are the $r$-tuples $b_1b_2\cdots b_r$
with $b_i\in \{0,1,\cdots,m_i-1\}$, and let two vertices be adjacent
if the corresponding tuples differ in precisely one place. Such a
graph is called a $Hamming~graph$. Clearly, a graph $G$ is a Hamming
graph if and only if it can be written in the form $G=K_{m_1}\square
K_{m_2}\cdots\square K_{m_r}$ for some $r\geq 1$, where $m_i\geq 2$
for all $i$. So we call $G$ a Hamming graph with $r$ factors. A
Hamming graph is a $hypercube$ (or $r$-$cube$) \cite{harary},
denoted by $Q_r$, if and only $m_i=2$ for all $i$. The concept of
Hamming graph is useful in communication networks \cite{Imrich}.

\begin{cor}\label{thm142}\cite{Li-Sun5}
If $G$ is a Hamming graph with $r$ factors, then $rc(G)=r.$ In
particular, $rc(Q_r)=r$.
\end{cor}

The $composition~(lexicographic~product)$ of two graphs $G$ and $H$
is the simple graph $G[H]$ with vertex set $V(G)\times V(H)$ in
which $(u,v)$ is adjacent to $(u',v')$ if and only if either $uu'
\in E(G)$ or $u=u'$ and $vv' \in E(H)$. By definition, $G[H]$ can be
obtained from $G$ by substituting a copy $H_v$ of $H$ for every
vertex $v$ of $G$ and by joining all vertices of $H_v$ with all
vertices of $H_u$ if $uv\in E(G)$. Note that $G[H]$ is connected if
and only if $G$ is connected. By definition, it is easy to show: If
$G$ is complete, then $diam(G[H])=1$ if $H$ is complete (as now
$G[H]$ is complete), $diam(G[H])=2$ if $H$ is not complete; If $G$
is not complete, then $diam(G[H])=diam(G)$. Then we have the
following result.

\begin{thm}\label{thm143}\cite{Li-Sun5}
If $G$ and $H$ are two graphs and $G$ is connected, then we have\\
(1) if $H$ is complete, then $$rc(G[H])\leq rc(G).$$ In
particular, if $diam(G)=rc(G)$, then $rc(G[H])= rc(G).$\\
(2) if $H$ is not complete, then $$rc(G[H])\leq rc(G)+1.$$ In
particular, if $diam(G)=rc(G)$, then $diam(G[H])=2$ if $G$ is
complete and $rc(G)\leq diam(G[H])\leq rc(G)+1$ if $G$ is not
complete.
\end{thm}

In \cite{Li-Sun5}, Li and Sun also investigated other graph
operations, such as the union of graphs which we will not introduce
here.

\subsection{An upper bound for strong rainbow connection number}

The topic of rainbow connection number is fairly interesting and
recently a series papers have been published about it. The strong
rainbow connection number is also interesting, and by definition,
the investigation of it is more challenging than that of the rainbow
connection number. However, there are very few papers that have been
published about it. In \cite{Chartrand 1}, Chartrand, Johns, McKeon
and Zhang determined the precise strong rainbow connection numbers
for some special graph classes including trees, complete graphs,
wheels and complete bipartite (multipartite) graphs as shown in
Subsection 2.1. However, for a general graph $G$, it is almost
impossible to give the precise value for $src(G)$, so we aim to give
upper bounds for it according to some graph parameters. Li and Sun
\cite{Li-Sun6} derived a sharp upper bound for $src(G)$ according to
the number of edge-disjoint triangles (if exist) in a graph $G$, and
give a necessary and sufficient condition for the sharpness. We need
to introduce a new graph class.

Recall that a $block$ of a connected graph $G$ is a maximal
connected subgraph without a cut vertex. Thus, every block of $G$ is
either a maximal 2-connected subgraph or a bridge. We now introduce
a new graph class. For a connected graph $G$, we say $G\in
\overline{\mathcal{G}}_t$, if it satisfies the following conditions:

\textbf{$C_1$}. \ Each block of $G$ is a bridge or a triangle;

\textbf{$C_2$}. \ $G$ contains exactly $t$ triangles;

\textbf{$C_3$}. \ Each triangle contains at least one vertex of
degree two in $G$.

\begin{figure}[!hbpt]
\begin{center}
\includegraphics[scale=1.000000]{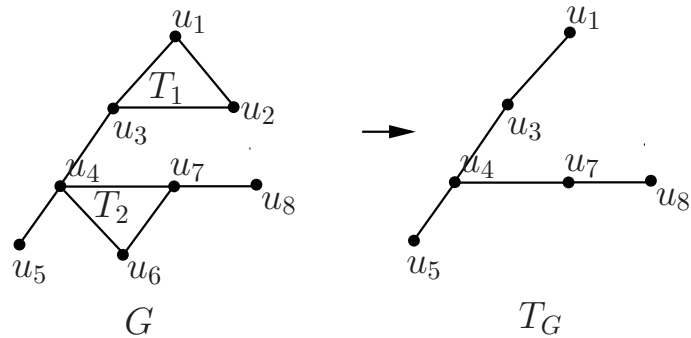}
\end{center}
\caption{An example of $G\in \overline{\mathcal{G}}_t$ with
$t=2$.}\label{figure4}
\end{figure}

By the definition, each graph $G\in \overline{\mathcal{G}}_t$ is
formed by (edge-disjoint) triangles and paths (may be trivial),
these triangles and paths fit together in a treelike structure, and
$G$ contains no cycles but the $t$ (edge-disjoint) triangles. For
example, see Figure \ref{figure4},
here $t=2$, and $u_1$, $u_2$ and $u_6$ are vertices of degree 2 in
$G$. If a tree is obtained from a graph $G\in
\overline{\mathcal{G}}_t$ by deleting one vertex of degree 2 from
each triangle, then we call this tree is a $D_2$-$tree$ of $G$,
denoted by $T_{G}$. For example, in Figure \ref{figure4}, $T_{G}$ is
a $D_2$-tree of $G$. Clearly, the $D_2$-tree is not unique, since in
this example, we can obtain another $D_2$-tree by deleting the
vertex $u_1$ instead of $u_2$. On the other hand, we can say that
any element of $\overline{\mathcal{G}}_t$ can be obtained from a
tree by adding $t$ new vertices of degree 2. It is easy to show that
the number of edges of $T_{G}$ is $m-2t$ where $m$ is the number of
edges of $G$.

They derived the following result.
\begin{thm}\label{thm144}\cite{Li-Sun6}
If $G$ is a graph with $m$ edges and $t$ edge-disjoint triangles,
then
$$src(G)\leq m-2t,$$ the equality holds if and only if $G\in
\overline{\mathcal{G}}_t$.
\end{thm}

In \cite{Ahadi}, Ahadi and Dehghan also derived an upper bound for
strongly regular graph: if $G$ is an $SRG(n,r,\lambda,\mu)$, then
$src(G)\leq \lceil
(e(4{\mu}r-4{\mu}{\lambda}-6{\mu}+1))^{\frac{1}{\mu}}\rceil$.
However, we do not know whether it is sharp.

Unlike rainbow connection number, which is a monotone graph property
(adding edges never increases the rainbow connection number), this
is not the case for the strong rainbow connection number (see Figure
\ref{figure2} for an example). The investigation of strong rainbow
connection number is much harder than that of rainbow connection
number. Chakraborty, Fischer, Matsliah and Yuster gave the following
conjecture.
\begin{con}\label{thm152}\cite{Chakraborty}
If $G$ is a connected graph with minimum degree at least $\epsilon
n$, then it has a bounded strong rainbow connection number.
\end{con}

\section{Rainbow $k$-connectivity}

In this section, we survey the results on rainbow $k$-connectivity.
By its definition, we know that it is difficult to derive exact
value or a nice upper bound of the rainbow $k$-connectivity for a
general graph. Chartrand, Johns, McKeon and Zhang \cite{Chartrand 2}
did some basic research on the rainbow $k$-connectivity of two
special graph classes. They first studied the rainbow
$k$-connectivity of the complete graph $K_n$ for various pairs $k,
n$ of integers, and derived the following result:
\begin{thm}\label{thm201}\cite{Chartrand 2} For every integer $k\geq
2$, there exists an integer $f(k)$ such that if $n\geq f(k)$, then
$rc_k(K_n)=2$.
\end{thm}

They obtained an upper bound $(k+1)^2$ for $f(k)$, namely $f(k)\leq
(k+1)^2$. Li and Sun \cite{Li-Sun3} continued their investigation,
and the following result is derived:
\begin{thm}\label{thm202}\cite{Li-Sun3}
For every integer $k \geq 2$, there exists an integer $f(k)=
ck^{\frac{3}{2}}+C(k)$ where $c$ is a constant and
$C(k)=o(k^{\frac{3}{2}})$ such that if $n \geq f(k)$, then
$rc_k(K_n)=2$.
\end{thm}

From Theorem \ref{thm202}, we can obtain an upper bound
$ck^{\frac{3}{2}}+C(k)$ for $f(k)$, where $c$ is a constant and
$C(k)=o(k^{\frac{3}{2}})$, that is, they improved the upper bound of
$f(k)$ from $O(k^2)$ to $O(k^{\frac{3}{2}})$, a considerable
improvement. Dellamonica, Magnant and Martin \cite{Dellamonica} got
the best possible upper bound $2k$, which is linear in $k$ (see
Theorem \ref{thm203}). However, the proof of Theorem \ref{thm202} is
more structural or constructive, and informative. In the argument of
\cite{Dellamonica}, Dellamonica, Magnant and Martin put forward a
new concept, the $rainbow~(k,l)$-$connectivity$.

Given an edge-colored simple graph $G$, let $l \leq k$ be integers.
Suppose the edges of $G$ are $k$-colored. For $a,b\in V(G)$, denote
by $p(a,b)$ the maximum number of internally disjoint rainbow paths
of length $l$ having endpoints $a$ and $b$. The
$rainbow~(k,l)$-$connectivity$ of $G$ is the minimum $p(a,b)$ among
all distinct $a,b\in V(G)$. Note that this new defined
$rainbow~(k,l)$-$connectivity$ computes the number of internally
disjoint paths with the $same$ length $l$ (this is distinct from the
rainbow $k$-connectivity which, as mentioned above, computes the
number of colors); and by definition, for different edge-colorings,
the values of $rainbow~(k,l)$-$connectivity$ could be different.

From Theorem \ref{thm201}, we know that for any $r$, there exists an
explicit 2-coloring of $K_r$ in which the number of bi-chromatic
paths of length 2 between any pair of vertices is at least $\lfloor
\sqrt{r}-1\rfloor$. Using the above definition, it is a statement
about the rainbow (2,2)-connectivity of a given 2-coloring of the
edges of $K_r$. In \cite{Dellamonica}, Dellamonica, Magnant and
Martin greatly improve and generalize the above lower bound for
graphs of sufficiently large order by providing a different
constructive coloring. Their construction attains asymptotically the
maximum rainbow connectivity possible.

\begin{thm}\label{thm203}\cite{Dellamonica}
For any $k\geq 2$ and $r\geq r_0=r_0(k)$ there exists an explicit
$k$-coloring of the edges of $K_r$ having rainbow
$(k,2)$-connectivity $$(\frac{k-1}{k}-o(1))r.$$
\end{thm}

More generally, They considered the problem of finding longer
rainbow paths.

\begin{thm}\label{thm204}\cite{Dellamonica}
For any $3\leq l\leq k$, there exists $r_0=r_0(k)$ such that for
every $r\geq r_0$, there is an explicit $k$-coloring of the edges of
$K_r$ having rainbow $(k,2)$-connectivity
$$(1-o(1)){\frac{r}{l-1}}.$$
\end{thm}

This result is also asymptotically best possible, since any
collection of internally disjoint paths of length $l$ can contain at
most $\frac{r}{l-1}$ paths. Their proof employed a very recent
breakthrough due to Bourgain \cite{Bourgain, rao}, which consists of
a powerful explicit extractor. Roughly speaking, an (explicit)
extractor is a polynomial time algorithm used to convert some
special probability distributions into uniform distributions. See
\cite{Shaltiel} for a good but somewhat outdated survey on
extractors.)

Chartrand, Johns, McKeon and Zhang \cite{Chartrand 2} also
investigated the rainbow $k$-connectivity of $r$-regular complete
bipartite graphs for some pairs $k, r$ of integers with $2\leq k\leq
r$, and they showed:
\begin{thm}\label{thm205}\cite{Chartrand 2}
For every integer $k\geq 2$, there exists an integer $r$ such that
$rc_k(K_{r,r})=3$.
\end{thm}

However, they could not show a similar result for complete graphs,
and therefore they left an open question: For every integer $k\geq
2$, determine an integer (function) $g(k)$, for which
$rc_k(K_{r,r})=3$ for every integer $r\geq g(k)$, that is, the
rainbow $k$-connectivity of the complete bipartite graph $K_{r,r}$
is essentially 3. In \cite{Li-Sun4}, Li and Sun solved this question
using a similar but more complicated method to that of Theorem
\ref{thm205}, and they proved:

\begin{thm}\label{thm206}\cite{Li-Sun4}
For every integer $k\geq 2$, there exists an integer
$g(k)=2k\lceil\frac{k}{2}\rceil$ such that $rc_k(K_{r,r})=3$ for any
$r\geq g(k)$.
\end{thm}

More generally, we have the following question.
\begin{que}\label{que2}\cite{Chartrand 2}
Does the following hold? For each integer $k\geq 2$, there exists an
integer $h(k)$ such that for every two integers $s$ and $t$ with
$h(k) \leq s\leq t$, we have $rc_k(K_{s,t})= 3$.
\end{que}

Recently, He and Liang \cite{He} investigated the rainbow
$k$-connectivity in the setting of random graphs. They determined a
sharp threshold function for the property $rc_k(G(n,p))\leq d$ for
every fixed integer $d\geq 2$. This substantially generalizes a
result due to Caro, Lev, Roditty, Tuza and Yuster (see Theorem
\ref{thm135}). Their main result is as follows.
\begin{thm}\label{thm220}\cite{He}
Let $d\geq 2$ be a fixed integer and $k=k(n)\leq O(\log n)$. Then
$p=\frac{(\log n)^{1/d}}{n^{(d-1)/d}}$ is a sharp threshold function
for the property $rc_k(G(n,p))\leq d$.
\end{thm}

The key ingredient of their proof is the following result: With
probability at least $1-n^{-\Omega(1)}$, every two different
vertices of $G(n,\frac{C(\log n)^{1/d}}{n^{(d-1)/d}})$ are connected
by at least $2^{10d}{c_0}{\log n}$ internally disjoint paths of
length exactly $d$.

So far, results on the rainbow $k$-connectivity are just those for
two special graph classes, and a sharp threshold function for the
property $rc_k(G(n,p))\leq d$ where $d\geq 2$ is a fixed integer and
$k=k(n)\leq O(\log n)$. Clearly, there are a lot of things one can
do further on this concept. As mentioned above, it is difficult to
derive the precise value or a nice upper bound for $rc_k(G)$ of a
$\kappa$-connected graph $G$, where $2\leq k\leq \kappa$. So one may
consider the following problem.

\begin{prob}\label{prob14}
Derive a sharp upper bound for $rc_2(G)$, where $G$ is a
$\kappa$-connected graph with $\kappa \geq 2$. Does $rc_2(G)\leq
{\alpha}n$, where $0<\alpha <1$ is independent of $n$?
\end{prob}

\section{$k$-rainbow index}

The $k$-rainbow coloring as defined above, is another generalization
of the rainbow coloring. In \cite{Chartrand 4}, Chartrand, Okamoto
and Zhang did some basic research on this topic. There is a rather
simple upper bound for ${r\textsl{x}}_k(G)$ in terms of the order of
$G$, regardless the value of $k$.
\begin{pro}\label{thm207}\cite {Chartrand 4} Let $G$ be a nontrivial
connected graph of order $n\geq 3$. For each integer $k$ with $3\leq
k\leq n-1$, ${r\textsl{x}}_k(G)\leq n-1$ while ${r\textsl{x}}_n(G)=
n-1$.
\end{pro}

There is a class of graphs of order $n\geq 3$ whose $k$-rainbow
index attains the upper bound of Proposition \ref{thm207}, it is an
immediate consequence of Observation \ref{ob1}.
\begin{pro}\label{thm208}\cite {Chartrand 4} Let $T$ be a tree of
order $n\geq 3$. For each integer $k$ with $3\leq k\leq n$,
${r\textsl{x}}_k(T)= n-1$.
\end{pro}

There is also a rather obvious lower bound for the $k$-rainbow index
of a connected graph $G$ of order $n$, where $3\leq k\leq n$. The
$Steiner~distance~d(S)$ of a set $S$ of vertices in $G$ is the
minimum size of a tree in $G$ containing $S$. Such a tree is called
a $Steiner~S$-$tree$ or simply a $Steiner~tree$. The
$k$-$Steiner~diameter$, say $sdiam_k(G)$ of $G$ is the maximum
Steiner distance of $S$ among all sets $S$ with $k$ vertices in $G$.
Thus if $k=2$ and $S=\{u,v\}$, then $d(S)=d(u,v)$ and the 2-Steiner
diameter $sdiam_2(G)=diam(G)$. The $k$-Steiner diameter then
provides a lower bound for the $k$-rainbow index of $G$: for every
connected graph $G$ of order $n\geq 3$ and each integer $k$ with
$3\leq k\leq n$, ${r\textsl{x}}_k(G)\geq sdiam_k(G)\geq k-1$.

\begin{thm}\label{thm209}\cite {Chartrand 4} If $G$ is a unicyclic
graph of order $n\geq 3$ and girth $g\geq 3$, then
\[
 {r\textsl{x}}_k(G)=\left\{
   \begin{array}{ll}
     n-2 &\mbox {if~$k=3~and~g\geq 4$,}\\
     n-1 &\mbox {if~$g=3~or~4\leq k\leq n$.}
   \end{array}
   \right.
\]
\end{thm}

The investigation of ${r\textsl{x}}_k(G)$ for a general $k$ and a
general graph $G$ is rather difficult. So one may consider
${r\textsl{x}}_k(G)$ for a special graph class, or, for a general
graph and small $k$, such as $k=3$.
\begin{prob}\label{prob15}
Derive a sharp upper bound for ${r\textsl{x}}_3(G)$.
\end{prob}

There is also a generalization of the $k$-rainbow index, say
$(k,\ell)$-$rainbow~index$ ${r\textsl{x}}_{k,\ell}$, of $G$ which is
mentioned in \cite{Chartrand 4} and we will not introduce it here.

\section{Rainbow vertex-connection number}

The above several parameters are all defined on edge-colored graphs.
Here we will introduce a new graph parameter which is defined on
vertex-colored graphs. It is, as mentioned above, a vertex-version
of the rainbow connection number. Krivelevich and Yuster
\cite{Krivelevich} put forward this new concept and proved a theorem
analogous to Theorem \ref{thm115}.

\begin{thm}\label{thm210}\cite {Krivelevich} A connected graph
$G$ with $n$ vertices has $rvc(G)< \frac{11n}{\delta(G)}$.
\end{thm}

The argument of this theorem used the concept of $k$-$strong$
two-step dominating sets. They proved that if $H$ is a connected
graph with $n$ vertices and minimum degree $\delta$, then it
contains a $\frac{\delta}{2}$-strong two-step dominating set $S$
whose size is at most $\frac{2n}{\delta+2}$. Then they derived an
edge-coloring for $G$ according to its connected
$\frac{\delta}{2}$-strong two-step dominating set. And they showed
that, with positive probability, their coloring yields a rainbow
connected graph by the Lov\'{a}sz Local Lemma (see \cite{Alon}).

Motivated by the method of Theorem \ref{thm210}, Li and Shi derived
an improved result.
\begin{thm}\label{thm211}\cite {Li-Shi1}
A connected graph $G$ of order $n$ with minimum degree $\delta$ has
$rvc(G)\leq 3n/(\delta +1)+5$ for $\delta \geq \sqrt{n-1}-1, n\geq
290$, while $rvc(G)\leq 4n/(\delta +1)+5$ for $16\leq \delta \leq
\sqrt{n-1}-2$, $rvc(G)\leq 4n/(\delta +1)+C(\delta)$ for $6\leq
\delta \leq 16$ where
$C(\delta)=e^{\frac{3\log({\delta}^3+2{\delta}^2+3)-3(\log3-1)}{\delta-3}-2}$,
$rvc(G)\leq n/2-2$ for $\delta =5$, $rvc(G)\leq 3n/5-8/5$ for
$\delta =4$,$rvc(G)\leq 3n/4-2$ for $\delta =3$. Moreover, an
example shows that when when $\delta \geq \sqrt{n-1}-1$, and
$\delta=3,4,5$ the bounds are seen to be tight up to additive
factors.
\end{thm}

Motivated by the method of Theorem \ref{thm210}, Dong and Li
\cite{Dong-Li} also proved a theorem analogous to Theorem
\ref{thm155} for the rainbow vertex-connection version according to
the degree sum condition $\sigma_2$, which is stated as the
following theorem.

\begin{thm}\label{thm224}\cite{Dong-Li}
For a connected graph $G$ of order $n$, $rvc(G)\leq
8{\frac{n-2}{\sigma_2+2}}+10$ for $2\leq \sigma_2\leq 6$ and
$\sigma_2\geq 28$, while for $7\leq \sigma_2\leq 8$ and $16\leq
\sigma_2\leq 27$, $rvc(G)\leq \frac{10n-16}{\sigma_2+2}+10$, and for
$9\leq \sigma_2\leq 15$, $rvc(G)\leq
\frac{10n-16}{\sigma_2+2}+A(\sigma_2)$, where $A(\sigma_2)=63, 41,
27, 20, 16, 13, 11$, respectively.
\end{thm}

Note that by the definition of $\sigma_2$, we know $\sigma_2\geq
2\delta$, so we have $8{\frac{n-2}{\sigma_2+2}}+10\leq
\frac{4(n-2)}{\delta+1}+10$, and hence the bound of Theorem
\ref{thm224} is an improvement to that of Theorem \ref{thm211} in
the case $16\leq \delta \leq \sqrt{n-1}-2$.

\section{Algorithms and computational complexity}

At the end of \cite{Caro}, Caro, Lev, Roditty, Tuza and Yuster gave
two conjectures (see Conjecture 4.1 and Conjecture 4.2 in
\cite{Caro}) on the complexity of determining the rainbow connection
numbers of graphs. Chakraborty, Fischer, Matsliah and Yuster
\cite{Chakraborty} solved these two conjectures by the following
theorem.

\begin{thm}\label{thm212}\cite{Chakraborty} Given a graph $G$, deciding if $rc(G)=2$ is
NP-Complete. In particular, computing $rc(G)$ is NP-Hard.
\end{thm}

Chakraborty, Fischer, Matsliah and Yuster divided the proof of
Theorem \ref{thm212} into three steps: the first step is showing the
computational equivalence of the problem of $rainbow$ $connection$
$number$ $2$, that asks for a red-blue edge coloring in which all
vertex pairs have a rainbow path connecting them, to the problem of
$subset$ $rainbow$ $connection$ $number$ $2$, that asks for a
red-blue coloring in which every pair of vertices in a given subset
of pairs has a rainbow path connecting them. In the second step,
they reduced the problem of $extending$ $to$ $rainbow$ $connection$
$number$ $2$, asking whether a given partial red-blue coloring can
be completed to obtain a rainbow connected graph, to the problem of
subset rainbow connection number 2. Finally, the proof of Theorem
\ref{thm212} is completed by reducing 3-SAT to the problem of
extending to rainbow connection number 2.

Chakraborty, Fischer, Matsliah and Yuster \cite{Chakraborty} also
raised the following problem.

\begin{prob}\label{prob10}\cite{Chakraborty}
Suppose that we are given a graph $G$ for which we are told that
$rc(G)=2$. Can we rainbow-color it in polynomial time with $o(n)$
colors ?
\end{prob}

For the usual coloring problem, this version has been well studied.
It is known that if a graph is 3-colorable (in the usual sense),
then there is a polynomial time algorithm that colors it with
$\tilde{O}(n^{3/14})$ colors \cite{blum}.

Suppose we are given an edge coloring of the graph. Is it then
easier to verify whether the colored graph is rainbow connected?
Clearly, if the number of colors is constant then this problem
becomes easy. However, if the coloring is arbitrary (with an
unbounded number of colors), the problem becomes NP-Complete.

\begin{thm}\label{thm213} \cite{Chakraborty} The following problem
is NP-Complete: Given an edge-colored graph $G$, check whether the
given coloring makes $G$ rainbow connected.
\end{thm}

For the proof of Theorem \ref{thm213}, Chakraborty, Fischer,
Matsliah and Yuster first showed that the $s-t$ version of the
problem is NP-Complete. That is, given two vertices $s$ and $t$ of
an edge-colored graph, decide whether there is a rainbow path
connecting them. Then they reduced the problem of Theorem
\ref{thm213} from it.

More generally it has been shown in \cite{Le}, that for any fixed
$k\geq 2$, deciding if $rc(G)=k$ is $NP$-complete.

In \cite{Chakraborty}, Chakraborty, Fischer, Matsliah and Yuster
also derived some positive algorithmic results. Parts of the
following two results were shown in Theorems \ref{thm114} and
\ref{thm136}. They proved:

\begin{thm}\label{thm214}\cite{Chakraborty} For every $\epsilon >0$
there is a constant $C=C(\epsilon)$ such that if $G$ is a connected
graph with $n$ vertices and minimum degree at least ${\epsilon}n$,
then $rc(G)\leq C$. Furthermore, there is a polynomial time
algorithm that constructs a corresponding coloring for a fixed
$\epsilon$.
\end{thm}

As mentioned above, Theorem \ref{thm214} is based upon a modified
degree-form version of Szemer\'{e}di¡¯s Regularity Lemma that they
proved and that may be useful in other applications. From their
algorithm it is also not hard to find a probabilistic polynomial
time algorithm for finding this coloring with high probability
(using on the way the algorithmic version of the Regularity Lemma
from  \cite{Alon2} or \cite{Fischer}).

\begin{thm}\label{thm215}\cite{Chakraborty} If $G$ is an $n$-vertex
graph with diameter 2 and minimum degree at least $8{\log}n$, then
$rc(G)\leq 3$. Furthermore, such a coloring is given with high
probability by a uniformly random $3$-edge-coloring of the graph
$G$, and can also be found by a polynomial time deterministic
algorithm.
\end{thm}

As mentioned, He and Liang \cite{He} investigated the rainbow
$k$-connectivity in the setting of random graphs. They also
investigated rainbow $k$-connectivity from the algorithmic point of
view. The NP-hardness of determining $rc(G)$ was shown by
Chakraborty et al. as shown above. They showed that the problem
(even the search version) becomes easy in random graphs.

\begin{thm}\label{thm221}\cite{He} For any constant $\epsilon \in [0,1)$,
$p=n^{-{\epsilon}(1\pm o(1))}$ and $k\leq O(\log n)$, there is a
randomized polynomial time algorithm that, with probability
$1-o(1)$, makes $G(n,p)$ rainbow-$k$-connected using at most one
more than the optimal number of colors.
\end{thm}

Since almost all natural edge probability functions $p$ encountered
in various scenarios have such form, their result is quite strong.
Note that $G(n,n^{-\epsilon})$ is almost surely disconnected when
$\epsilon >1$ \cite{Erdos2}, which makes the problem become trivial.
Therefore they ignored these cases.

Recall that the values of rainbow~$(k,l)$-connectivity may be
different for distinct edge-colorings. In \cite{Dellamonica},
Dellamonica, Magnant and Martin derived that a random $k$-coloring
of a sufficiently large complete graph has asymptotically optimal
rainbow rainbow~$(k,l)$-connectivity (see Theorems \ref{thm203} and
\ref{thm204}). They obtained an $explicit$ edge-coloring. By
explicit, we mean that they gave a polynomial time algorithm to
compute such an edge-coloring.

Recently, the computational complexity of rainbow vertex-connection
numbers has been determined by Chen, Li and Shi \cite{Chen-Li-Shi}.

Motivated by the proofs of Theorems \ref{thm212} and \ref{thm213},
they derived two corresponding results to the rainbow
vertex-connection.

\begin{thm}\label{thm223}\cite{Chen-Li-Shi} Given a graph $G$, deciding if $rvc(G)=2$ is
NP-Complete. In particular, computing $rvc(G)$ is NP-Hard.
\end{thm}

\begin{thm}\label{thm222}\cite{Chen-Li-Shi}
The following problem is NP-Complete: given a vertex-colored graph
$G$, check whether the given coloring makes $G$ rainbow
vertex-connected.
\end{thm}

\end{document}